\documentclass[12pt]{article} 
 
\usepackage{latexsym} 
\usepackage{ntheorem} 
\usepackage[latin1]{inputenc} 

\usepackage{amsfonts,amsmath,amssymb} 

\usepackage{a4wide}

\author{Vincent Bosser (Caen)\footnote{Supported by the contract ANR ``HAMOT", BLAN-0115-01.} \\ 
Andrea Surroca (Basle)\footnote{Supported by an Ambizione fund PZ00P2\_121962 of the Swiss National Science Foundation and the Marie Curie IEF 025499 of the European Community.}}

\title{\Large Elliptic logarithms,  diophantine approximation \\
and the Birch  and Swinnerton-Dyer  conjecture} 
   
\date{\today} 
 
\newtheorem{thm}{\textbf{Theorem}}[section]
\newtheorem{conj}[thm]{\textbf{Conjecture}}
\newtheorem{lemma}[thm]{\textbf{Lemma}} 
\newtheorem{prop}[thm]{\textbf{Proposition}}

{\theorembodyfont{\upshape}
\newtheorem{remark}[thm]{Remark}}

\newcommand\rat{\mathbf{Q}} 
\newcommand\Q{\mathbf{Q}} 
\newcommand\C{\mathbf{C}} 
\newcommand\Z{\mathbf{Z}} 
\newcommand\enteros{\mathbf{Z}} 
\newcommand\PP{\mathbf{P}} 
\newcommand\R{\mathbf{R}} 
 
\newcommand\Kbarre{\overline{K}}

\newcommand\card{\mathrm{card}}
\newcommand\rad{\mathrm{rad_F}}

\newcommand\ord{\mathrm{ord}}
\newcommand\rk{\mathrm{rk}}

\newcommand\Reg{\mathrm{Reg}}

\newcommand\Gal{\mathrm{Gal}} 
 
\newcommand\vp{\mathfrak{p}} 
\newcommand\vq{\mathfrak{q}} 
\newcommand\va{\mathfrak{a}} 

\newcommand\E{\mathcal{E}} 
 
\newcommand\Diff{\Omega_{\mathcal{E}/O_{K}}^{1}}

\newcommand\Fu{{K_{0}}}
\newcommand\Fd{F} 
\newcommand\Ft{K} 
  
\newcommand\espproj{\mathbf{P}} 
 
\newcommand\esp{\hspace{0,2cm}} 
\newcommand\Belyi{Bely\u\i}

\input cyracc.def \font\tencyr=wncyr10 \def\russe{\tencyr\cyracc} 
\def\Sha{\text{\russe{Sh}}} 
 
\begin{document} 
 
\bibliographystyle{alpha}

\maketitle 
 
\begin{quote} 
\textbf{Abstract.} 
{\small   
Most, if not all, unconditional results towards the $abc$-conjecture rely ultimately on classical Baker's method. 
In this article, we turn our attention to its {\it elliptic} analogue.
Using the elliptic Baker's method, we have recently obtained a
new upper bound for the height of the $S$-integral points on an elliptic curve.
This bound depends on some parameters related to the Mordell-Weil group of the curve.
We deduce here a  bound relying on the conjecture of Birch and Swinnerton-Dyer,
involving classical, more manageable quantities.
We then study which $abc$-type inequality over number fields could be derived from this elliptic approach.
} 

\end{quote}
 
\begin{quote} 
\textbf{Keywords:} 
{\small  integral points on elliptic curves, quantitative Siegel's theorem, elliptic logarithms, Birch and Swinnerton-Dyer conjecture, $abc$-conjecture.} 
\end{quote}

\begin{quote} 
\textbf{Mathematical subject classification:} 
{\small  Primary: 11G50; Secondary: 11G05, 11J86, 14G05, 11G40.} 
\end{quote}
 

\section{Introduction} 

The $abc$-conjecture of D. W. Masser  and J. Oesterl\'e (Conjecture \ref{abc} below) is one of the most important unsolved problems in Diophantine
analysis. It is well known that it is connected with several problems in number theory. If true, it would imply strong or quantitative versions
of important theorems. 
Indeed, let us recall the classical Siegel's theorem: for an affine curve of genus $\geq 1$ or of genus $0$ having at least three points at infinity,
the set of integral points is finite. This theorem was later extended to $S$-integral points by Mahler. 
For curves of genus $\geq 2$, Siegel's theorem is superseded by Faltings' theorem, which asserts that the set of rational points on an algebraic
curve of genus greater than 2 is finite. These results are qualitative statements in general, that is, there is no known proof providing an upper
bound for the height of the points, \textit{i.e.} a ``quantitative" result, which would allow to find these points. The only known quantitative
results on rational points concern integral points 
and only some particular cases, \textit{e.g.} the case of curves of genus $0$, $1$, or the case of curves which are Galois coverings of
$\espproj^1 \setminus \{0,1, \infty \}$.
All of them come from classical Baker's method using lower bounds for linear forms in logarithms.
We refer the reader to \cite{gyory-panorama}, \cite{bilu-panorama} and the references therein for an overview of known results.

As noticed by Elkies \cite{elkies-abc}, the $abc$-conjecture over number fields would imply a quantitative version on Faltings' theorem.
As shown by the second author \cite{parma}, also a quantitative Siegel's theorem would follow, with explicit dependence on the set of places $S$.
Unfortunately only weak results are known towards this conjecture yet, and they are insufficient to yield a quantitative Siegel's theorem
for new classes of curves.

On the other hand, it is worth noting that Moret-Bailly \cite{moret-bailly} showed that, conversely, a uniform and effective version of Falting's theorem for the curve $y^2+y=x^5$ would imply $abc$.  As shown in \cite{bornes}, any bound for the height of the more restrictive set of the $S$-integral points
on a fixed curve, explicit in the set $S$ and in the degree and the discriminant of the number field considered, would suffice to imply a result
towards the $abc$-conjecture over this number field.
Using such a bound given by a quantitative Siegel's theorem due to Bilu \cite{bilu.compositio}, the second author obtained in her thesis
(see \cite{bornes}) the first result towards the $abc$-conjecture over an arbitrary number field.
Afterwards, K. Gyory and K. Yu (\cite{gyory-yu}, \cite{gyory.acta-arith2008}) gave completely explicit $abc$-type results using
bounds for the height of solutions of $S$-unit equations. When the number field is $\Q$ better inequalities were known, see \cite{stewart-tijdeman} for the first result obtained and 
\cite{stewart-yu.1991}, \cite{stewart-yu.2001} for later improvements.
Roughly speaking, all these results differ from the conjecture from an exponential.\footnote{
In another direction, A. Baker \cite{baker.abc} and P. Philippon \cite{philippon} suggested some 
conjectures on linear forms in logarithms which would imply a weak version 
of $abc$ (where $1+\epsilon$ is replaced by some other constant).
These conjectures involve simultaneously several places (archimedean and non-archimedean) but these kinds of results are now far away from being proved.
}

All the quantitative Siegel's theorems and the $abc$-type results mentioned above depend ultimately on lower bounds for
linear forms in usual logarithms, complex as well as $p$-adic. In this paper, we turn our attention to the {\it elliptic} analogue of
Baker's method. In fact, in order to get a quantitative Siegel's theorem in the case of an elliptic curve,
it seems to be more natural to take advantage of the group law
and thus to use linear forms in \textit{elliptic} logarithms. Following this approach,
we have recently obtained in \cite{S-entiers-Bosser-Surroca} a new upper bound for the height of the
$S$-integral points of an elliptic curve $E$ defined over a number field $K$,
using the explicit lower bounds for linear forms in elliptic logarithms of S. David \cite{david} in the archimedean case, and of 
N. Hirata \cite{noriko} in the ultrametric one. However, the method leads to a bound
which is not quite effective since it depends on several parameters depending on
the Mordell-Weil group $E(K)$ of the curve (see Theorem~\ref{borne-hauteur} below).

This raises the question of knowing under which conditions one can get an explicit, more manageable, upper bound in terms of the set $S$ and
the number field $K$ using the elliptic Baker's method, and which kind of result towards the $abc$-conjecture can be obtained in this way. 
This paper gives an answer to these two questions. 
We will see that, following Manin's conditional algorithm \cite[Theorem 11.1]{manin} based on the assumption of the conjecture of B. J. Birch and H. P. F. Swinnerton-Dyer (Conjecture~\ref{bsd}, BSD-conjecture for short) and the classical Hasse-Weil conjecture,
we can derive from \cite{S-entiers-Bosser-Surroca} a quantitative Siegel's theorem whose bound is explicit in $S$, the degree and the discriminant of the number field
(Theorem~\ref{conjub}). 
Thus this paper highlights some connection between Baker's method, the BSD-conjecture and
the $abc$-conjecture.


\medskip 
 
The paper is organized as follows. For convenience to the reader, we have gathered in Section~\ref{notations} the notation which will be used
throughout the text. In Section~\ref{mordell-weil}, after recalling the BSD-conjecture, we state and prove
a conjectural upper bound for the height of $S$-integral points (Theorem~\ref{conjub}).
Finally, in Section \ref{appli-abc}, we
prove Theorem \ref{bsd-abc} and discuss the result.

\section{Notations}\label{notations}

Throughout the text, if $x$ is a non negative real number, we set $\log^+x=\max\{1,\log x\}$ (with the convention $\log^+0=1$).
  
If $K$ is a number field, we will denote by $O_K$ its ring of integers, by
$D_K$ the absolute value of its discriminant, and by $M_K$ the set of places
of $K$. The set of the archimedean places will be denoted by $M_K^{\infty}$
and the set of the ultrametric ones will be denoted by $M_K^0$. For each $v$
in $M_K$, we define an absolute value $|\cdot|_v$ on $K$ as follows.
If $v$ is archimedean, then $v$ corresponds to an embedding 
$\sigma : K \hookrightarrow \C$ or its complex conjugate (we identify the place $v$ with
the embedding $\sigma$), and we set $|x|_v=|x|_{\sigma}:=|\sigma(x)|$, where
$|\cdot |$ is the usual absolute value on $\C$. If $v$ is ultrametric, then
$v$ corresponds to a non zero prime ideal $\mathfrak p$ of $O_K$ (we
will identify $v$ and ${\mathfrak p}$),
and we take for
$|\cdot|_v=|\cdot|_{\mathfrak p}$ the absolute value on $K$ normalized
by $|p|_v=p^{-1}$, where $p$ is the prime number such that
${\mathfrak p}\mid p$.
We denote by $K_v$ the completion of $K$ at $v$ and use
again the notation $|\cdot|_v$ for the unique extension of $|\cdot|_v$ to
$K_v$. If $v$ is an ultrametric place associated to the prime ideal
${\mathfrak p}$, we denote by $e_{\mathfrak p}$ the ramification index of
${\mathfrak p}$ over $p$, by $f_{\mathfrak p}$ the residue class degree,
and by $\ord_{\mathfrak p} : K _{\mathfrak{p}}^*\rightarrow\Z$ the valuation
normalized by $\ord_{\mathfrak p}(p)=e_{\mathfrak p}$
(hence $\ord_{\mathfrak p}(x)=-e_{\mathfrak p}\log_p|x|_{\mathfrak p}$ for
all $x$ in $K_{\mathfrak{p}}^*$).
 
If $S$ is a finite subset of $M_K^0$, we denote by 
$$O_{K,S} = \{x \in K; \forall v \notin S\cup M_K^{\infty}, 
|x|_{v} \leq 1  \}$$ 
the ring of $S$-integers of $K$, and we set 
$$\Sigma_S =\sum_{{\mathfrak p}\in S}\log N_{\Ft/\rat}({\mathfrak p}).$$ 
Note that with our notation, the set $S$ contains only non-archimedean places of $K$.
 
Throughout the text, we denote by $h$ the absolute logarithmic Weil height
on the projective space $\PP^n(\overline{\rat})$, and we denote by
$h_K:=[K:\rat]h$ the relative height on $\PP^n(K)$. Thus,
if $(\alpha_{0}: \ldots: \alpha_{n}) \in \mathbf{P}^{n}(K)$, we have:
\begin{equation}\label{h.weil} 
h(\alpha_{0}: \ldots : \alpha_{n}) = \frac{1}{[K:\rat]} 
\sum_{v \in M_{K}} [K_v:\Q_v] \log \max \{|\alpha_{0}|_{v}, \ldots, |\alpha_{n}|_{v}\}. 
\end{equation} 

For every $(\alpha_{1}: \alpha_{2}: \alpha_{3}) \in \mathbf{P}^{2}(K)$, we further define
$${\mathrm{rad_K}}(\alpha_{1}: \alpha_{2}: \alpha_{3}) = \Sigma_S,$$
where $S = \{\vp \in M_{K}^0; \card\{\ord_{\vp}(\alpha_{1}), \ord_{\vp}(\alpha_{2}), \ord_{\vp}(\alpha_{3})\}\geq 2 \}$.
 
\medskip 

Let $E$ be an elliptic curve defined over a number field $K$.
The Mordell-Weil group  $E(K)$ of $K$-rational points of $E$ is a finitely generated group:
$$E(K) \simeq E(K)_{tors} \oplus \enteros^{\rk(E(K))}.$$ 
We will often simply write $r = \rk(E(K))$ for its rank, and we will
denote by $(Q_1,\ldots,Q_r)$ a basis of its free part.

We  further denote by $\hat{h}:E(\overline{K})\rightarrow\R$ the N\'eron-Tate height on $E$.
The ``N\'eron-Tate pairing'' $<\, ,\, >$ is defined by $<P,Q> = \frac{1}{2}(\hat{h}(P+Q) - \hat{h}(P) - \hat{h}(Q))$.
The regulator $\Reg(E/K)$ of $E/K$ is the determinant of the matrix $\mathcal{H}=(<Q_{i}, Q_{j}>)_{1\leq i, j\leq r}$ of the N\'eron-Tate
pairing with respect to the chosen basis $(Q_1,\ldots,Q_r)$, that is 
$$\Reg(E/K) =  \det(\mathcal{H}).$$
  
Suppose now that the elliptic curve $E$ is embedded in $\PP^2$ and given by a Weierstrass equation
\begin{equation}\label{weierstrass} 
y^{2} = x^{3} + Ax + B 
\end{equation} 
with $A,B\in O_K$. Let us denote by $O=(0:1:0)$ the zero element of $E(K)$.
If $Q\ne O$ is a point of $E$, we will denote its affine coordinates
(in the above Weierstrass model) as usual by $(x(Q),y(Q))$.
For $Q$ in $E(\overline{\Ft})$ we define $h_x(Q):=h(1:x(Q))$ if $Q\not=O$ and $h_x(O):=0$. 
Finally, we denote by $E(O_{K,S})$ the set of $S$-integral points of $E(K)$
with respect to the $x$-coordinate, that is
$$E(O_{K,S}) = \{ Q \in E(K)\setminus\{O\}; 
x(Q) \in O_{K,S}  \} \cup \{O\}.$$

In fact, in all what follows it will be crucial to distinguish the field of definition of the elliptic curve from
the field of rationality of the points we will consider. More precisely, we will fix a number field $K_0$
and an elliptic curve $E$ defined over $K_0$, and we will consider points in $E(K)$, where $K$ is a finite extension of $K_0$
(that we will think of as varying). In the estimates that will occur we will neither be interested in the dependence on $E/K_0$ nor
try to explicit it, and we will thus consider as a ``constant'' any quantity depending on $E/K_0$.
This convention about constants will apply in particular to the various implicit constants involved in the symbols $\ll$ appearing in
the text.


\section{Conditional upper bound for the height of $S$-integral points }\label{mordell-weil} 
 
In this section, we fix a number field $K_0$ and we consider an elliptic curve $E$ defined over $K_0$ given
by a Weierstrass equation (\ref{weierstrass}), where $A,B \in O_{K_0}$.
Let $\Ft$ be a finite extension of $\Fu$ and $S\subset M_{\Ft}^0$ be a finite set of places of $K$.
We denote by $r$ the rank of the Mordell-Weil group $E(K)$, 
by $(Q_1,\ldots, Q_r)$ a system of generators of its free part, and we 
define the real number $V$ by
$$\log V:=\max\{\hat{h}(Q_i); 1\le i\le r\}.$$
We further set 
$$d:=[\Ft:\Q].$$

\smallskip 
 
In \cite{S-entiers-Bosser-Surroca}, we have obtained the following result. 

\begin{thm}\label{borne-hauteur}  In the above set up, let $Q$ be a point in $E(O_{\Ft,S})$. Then there exist positive effectively computable real numbers
$\gamma_0, \gamma_1$ and $\gamma_2$ depending only on $A$ and $B$ (that is, on the curve $E/\Fu$), such that,
if $r = 0$, then $h_x(Q) \le \gamma_0$, and, if $r > 0$, then
\begin{equation}\label{bound} 
h_x(Q) \le C_{E,\Ft}e^{(8r^2+\gamma_1dr)\Sigma_S},
\end{equation} 
where 
\begin{align}\label{bound2}
C_{E,\Ft} & = \gamma_2^{r^2}r^{2r^2} d^{9r+15}(\log^+d)^{r+6} 
(\log^+\log V)^{r+7}  (\log^+\log^+\log V)^2\, \prod_{i=1}^{r}\max\{1,\hat h(Q_i)\} \cr 
& \times \log^+\!(\Reg(E/K)^{-1}) 
(\log^+\!\log(\Reg(E/K)^{-1}))^2 
(\log^+\!\log^+\!\log(\Reg(E/K)^{-1})).
\end{align}  
\end{thm}

The aim of this section is to deduce from Theorem~\ref{borne-hauteur}
an upper bound for the height of the $S$-integral points of $E(K)$, depending only on $E/\Fu$ and on
the parameters $\Sigma_S$, $d$ and $D_K$.
Such a bound will be obtained assuming the Hasse-Weil conjecture and BSD-conjecture. The approach relies on Manin's algorithm \cite{manin}. See also  Masser's book \cite[Appendix IV, p. 140]{davidmasser}, where the association of Manin's algorithm with linear forms in elliptic logarithms appears for the first time to get an effective version of Siegel's theorem. The precise statement we prove here is given in the next section.

\subsection{Statement of the result}

In order to state the conjectural quantitative Siegel's theorem we obtain, we need to introduce the two conjectures we will use.
Denote by $L(E/K, s)$ the $L$-series (or $\zeta$-function) of $E/K$ at $s$, which is an analytic function for all $s$ satisfying $\Re(s) > \frac{3}{2}$. 
Let $\mathcal{F}_{E/K}$ denote the conductor of $E$ over $K$.
Following \cite{milne}, define the normalized $L$-function as 
$$\Lambda(E/K,s) = N_{K/\rat}(\mathcal{F}_{E/K})^{s/2} \cdot D_K^s\cdot ((2\pi)^{-s} \cdot \Gamma (s))^{[K:\rat]} \cdot L(E/K, s).$$ 
We then have the classical conjecture. 
 
\smallskip 
 
\begin{conj} [Hasse-Weil] \label{funct-eq} 
Let $E/K$ be an elliptic curve defined over a number field. The $L$-series of $E/K$ has an analytic continuation of finite order to the entire complex plane and satisfies the functional equation  
$$\Lambda(E/K, 2-s) = \varepsilon \Lambda (E/K, s), \esp \textrm{for some }\esp \varepsilon = \pm 1.$$ 
\end{conj} 
 
This conjecture is true for abelian varieties with complex multiplication (\cite{shimura-taniyama}), for elliptic curves over $\rat$ (\cite{wiles1995}) and in some special cases, this conjecture is also true for modular abelian varieties (\cite{shimura.automorphic.book}). 
 
\smallskip 
 
Denote $\Sha (E/K)= \ker (H^1(\Gal(\Kbarre/K), E_K) \rightarrow \prod_v H^1(\Gal(\Kbarre_v/K_v), E_{K_v}))$ the Tate-Shafarevich group, which is conjectured to be a finite group. (See \cite{rubin1987} and \cite{kolyvagin1988} for the first examples of elliptic curves for which it was proved that the Tate-Shafarevich group is finite.) Let $F(x,y) = 0$ be a Weierstrass equation for $E$. Denote $F_{y}$ the partial derivative of $F$ with respect to $y$. Then the invariant differential of the Weierstrass equation, $\omega = \frac{dx}{F_{y}}$,  
is holomorphic and non vanishing. 
Let $\E$ denote the N\'eron model of $E$ over $O_K$ and let $\Diff$ be the invertible sheaf of the differential 1-forms on $\E$. The module $H^{0}(\E, \Diff)$ of global invariant differentials on  $\E$ is  a projective $O_{K}$-module of rank 1 and can be written as 
$$H^{0}(\E, \Diff) = \omega \va,$$ 
where $\va$ is a fractional ideal of $K$ (depending on $\omega$). 
 
To every place $v$ of $K$, we will associate a local number  $c_v$.  For $v$ a finite place of $K$, let $E^0(K_v)$ be the subgroup of $K_v$-rational points which reduces to the identity component of the N\'eron model $\mathcal{E}$. Let $\mu_v$ be an additive Haar measure on $K_v$ such that $\mu_v(O_{K_v}) = 1$ if $v$ is finite, $\mu_v$ is the Lebesgue measure if $v$ is a real archimedean place and twice the Lebesgue measure if $v$ is complex. Define, for an {\it archimedean} place $v$, the {\it local period}
$$c_v = \int_{E(K_{v})}|\omega| \mu_v.$$
 Define, for a {\it finite} place $v$ of $K$, $c_v = \card(E(K_{v})/E^{0}(K_{v}))$ and   
$$c_{\infty}(E/K) = \prod_{v \in M_K^{\infty}}c_v \cdot N_{K/\rat}(\va),$$ 
which is independent of the choice of the differential $\omega$. 
 
The Birch and Swinnerton-Dyer conjecture can be stated as follows \cite{birch-swinnerton-dyer} (see also \cite{gross}).
 
\begin{conj}[Birch and Swinnerton-Dyer]\label{bsd} 
Let $E/K$ be an elliptic curve defined over a number field.
\begin{enumerate} 
\item The $L$-series $L(E/K,s)$ has an analytic continuation to the entire complex plane. 
\item $\ord_{s=1}L(E/K,s) = \rk (E(K))$. 
\item The leading coefficient $L^{\star}(E/K, 1) = \lim_{s\to1}\frac{L(E/K,s)}{(s-1)^{\rk(E(K))}}$ in the Taylor expansion of $L(E/K,s)$ at $s=1$ satisfies  
 \begin{equation}\label{formulaBSD}
 L^{\star}(E/K, 1) =   | \Sha (E/K)| \cdot \Reg(E/K) \cdot |E(K)_{tors}|^{-2}   \cdot c_{\infty}(E/K) \cdot \prod _{v \in M_K^0}c_v  \cdot D_K^{-1/2}. 
 \end{equation}
\end{enumerate} 
\end{conj}

\smallskip

We can now state the conjectural quantitative Siegel's theorem that we obtain.

\begin{thm}\label{conjub}
Let $\Fu$ be a number field, and let $E$ be an elliptic curve given by a Weierstrass equation $y^2=x^3+Ax+B$ with
$A, B\in O_{\Fu}$.  Let $K/\Fu$ be a finite extension, $S$ a finite set of
finite places of $K$, and denote by $d$ the degree $[K:\rat]$.
 
Suppose that the $L$-series of $E/K$ satisfies a  functional equation (Conjecture \ref{funct-eq}) 
and that the Birch and Swinnerton-Dyer Conjecture (Conjecture \ref{bsd}) holds for $E/K$.

Then, there exist positive numbers $\alpha_1$ and  $\alpha_2$ (depending on $E/K_0$ only) such that,  for every point $Q$ in $E(O_{K,S})$, we have
$$h_x(Q) \leq \exp\{\alpha_1^d +\alpha_2 \,d^6 (\log^+ D_K)^2\, [\Sigma_S + \log (d \log^+ D_K)]\}.$$
\end{thm}

The rest of Section~\ref{mordell-weil} is devoted to the proof of this theorem.
To deduce Theorem~\ref{conjub} from  Theorem~\ref{borne-hauteur}, we need to bound from above in terms of $d$, $D_K$ and $\Sigma_S$ the following parameters :
the rank $r$, the product $\prod_{i=1}^{r}\max\{1,\hat{h}(Q_{i})\}$, the greatest height $\log V$ and
the inverse of the regulator $\Reg(E/K)^{-1}$.
In Section~\ref{ub-rank}, we first bound from above the rank $r$. Then, in Section~\ref{ub-generators}, we bound from above
the three remaining quantities. The bounds for $\prod_{i=1}^{r}\max\{1,\hat{h}(Q_{i})\}$
and $\log V$ will rely on the BSD-conjecture.
Finally, we prove Theorem~\ref{conjub} in Section~\ref{conj-u_b}.

\subsection{An upper bound for the rank of the Mordell-Weil group}\label{ub-rank}
 
Explicit computations for the Weak Mordell-Weil theorem give a bound for the rank of the Mordell-Weil group in terms of the discriminant
of the number field.\footnote{Remark that, contrary to our situation, most often the interest in bounding the rank lies in the dependence
on the conductor of the elliptic curve, and the dependence on the number field is not considered. For example, under
Conjectures~\ref{funct-eq} and \ref{bsd} and the generalized Riemann hypothesis for $L(E/\rat)$, one would obtain
$\rk(E(K)) \ll \frac{\log \mathcal{F}_{E/K}}{\log \log \mathcal{F}_{E/K}}$, where the implied constant in $\ll$ depends on $K$.
(See \cite{mestre1986Compo}.)}
The following is Theorem 1 of \cite{ooe-top},  slightly corrected by G. R\'emond (\cite[Proposition 5.1]{gaelremondpreprint08}), 
for the special case where the abelian variety is an elliptic curve. We denote by $\mathcal{F}^0_{E/\Ft}$ the radical of the conductor of $E$ over $\Ft$, that is, the product of the prime ideals of $O_K$ where $E$ has bad reduction.
 
\begin{lemma}[Ooe-Top, R\'emond]\label{rank}   
There exist real numbers $\kappa_{1}, \kappa_{2}$ and $\kappa_{3}$ depending only on the degree $d = [\Ft:\rat]$ such that 
$$\rk E(\Ft) \leq \kappa_{1} \log N_{\Ft/\rat} \mathcal{F}^0_{E/\Ft}   +\kappa_{2} \log D_{\Ft} + \kappa_{3},$$ 
and one may take 
$\kappa_{2} = \frac{2^7}{\log 2}d, \, \, \,\kappa_{1} = 4d\kappa_{2} \, \, \,  {\textrm and} \, \,  \, \kappa_{3} = \kappa_{2} (\log 16 \cdot d^2 - 1 )$.
 \end{lemma}

\begin{lemma}\label{conductor} 
The conductors of the elliptic curve over $\Ft$ and over $\Fu$ satisfy 
$$\log N_{\Ft/\Q} \mathcal{F}_{E/\Ft} \leq 8 [\Ft:\Q] \log N_{\Fu/\rat} \mathcal{F}_{E/\Fu}.$$ 
\end{lemma} 
 
\noindent {\it Proof. } 
We will use the upper bounds for the exponents of the conductor  $\mathcal{F}_{E/\Ft}$ given in \cite[Theorem~0.1]{LRS} and \cite[Theorem~6.2]{BK}. According to these results, if we write the conductor $\mathcal{F}_{E/\Ft}$ of $E/\Ft$ as 
$$\mathcal{F}_{E/\Ft} = \prod _{\vq} \vq^{\delta_{\vq}(E)},$$ 
then, for each prime ideal $\vq$ lying above the prime number $p$, we have 
$\delta_{\vq}(E) \leq 2$ if $p\ge 5$, 
$\delta_{\vq}(E)\leq 2+3e_{\vq}$ if $p=3$, and 
$\delta_{\vq}(E)\leq 2+6e_{\vq}$ if $p=2$. 
In all cases we thus have the bound $\delta_{\vq}(E) \leq 8e_{\vq}$. Moreover, it is well known that if $\delta_{\vq}(E)\not= 0$, then $\vq$ must lie above a prime ideal ${\mathfrak p}$ of $O_{\Fu}$ at which $E$ has bad reduction. Since the prime ideals $\vp$ of bad reduction are those dividing $\mathcal{F}_{E/\Fu}$, we obtain 
\begin{align*} 
 N_{\Ft/\rat}(\mathcal{F}_{E/\Ft}) & = \prod_{\vp | \mathcal{F}_{E/\Fu}} 
\prod_{\vq | \vp} N_{\Ft/\rat}(\vq)^{\delta_{\vq}(E)} 
\leq \prod_{\vp | \mathcal{F}_{E/\Fu}}\prod_{\vq | \vp} 
N_{\Fu/\rat}(\vp)^{8e_{\vq}f_{\vq/\vp}} \cr 
& \leq \prod_{\vp | \mathcal{F}_{E/\Fu}}N_{\Fu/\rat}(\vp)^{8[\Ft:\rat]} 
\leq N_{\Fu/\rat}(\mathcal{F}_{E/\Fu})^{8[\Ft:\rat]}. 
\end{align*} 
\hfill $\Box$ 
 
Lemmas~\ref{rank} and \ref{conductor}, together with the fact that
$N_{\Ft/\rat} \mathcal{F}^0_{E/\Ft} \leq N_{\Ft/\rat} \mathcal{F}_{E/\Ft}$, lead to the following bound for the rank. 

\begin{lemma}\label{rank-discriminant}
The rank $r$ of the Mordell-Weil group $E(K)$ satisfies
$$r \ll d^3 (\log^+ D_K),$$
where the implicit constant depends at most on $E/K_0$.
\end{lemma}

\medskip

\subsection{On the generators of the Mordell-Weil group}\label{ub-generators}

We give here  upper bounds for $\log V$, $\prod_{i=1}^{r}\max\{1, \hat{h}(Q_{i})\}$ and $\Reg(E/K)^{-1}$.
For this purpose, we follow the approach of Yu. Manin \cite{manin}. We argue in two steps. In Section~\ref{unconditional-bounds}, we obtain first unconditional upper bounds, but involving the regulator $\Reg(E/K)$.
The main ingredients used are a result on the geometry of numbers as well as a lower bound for the N\'eron-Tate height of non-torsion points
due to Masser. Then, in Section~\ref{conditional-bounds}, we bound from above the regulator $\Reg(E/K)$ using the BSD-conjecture.

\subsubsection{An upper bound for $\log V$, $\prod_{i=1}^{r}\max\{1, \hat{h}(Q_{i})\}$ and $\Reg(E/K)^{-1}$}\label{unconditional-bounds}

\medskip
\noindent {\bf Geometry of numbers.}  
Recall that the N\'eron-Tate height on $E$ extends to a positive definite quadratic form on $E(K) \otimes_{\mathbf{Z}} \mathbf{R}$. The N\'eron-Tate pairing  gives $E(K) \otimes_{\mathbf{Z}}\mathbf{R}\simeq \mathbf{R}^r$ the structure of an Euclidean space. The free part of the Mordell-Weil group, $\Lambda:=E(K)/E(K)_{tors},$ sits as a lattice in this vector space. The regulator of $E/K$ is the square of the volume of a fundamental domain for the lattice. Thus we have
$$\Reg(E/K) = (\det(\Lambda))^2 = \det(\mathcal{H}),$$
where $\mathcal{H}=(<Q_{i}, Q_{j}>)_{1\leq i, j\leq r}$ is the matrix
of the N\'eron-Tate pairing with respect to the chosen basis $(Q_1,\ldots,Q_r)$. We begin by choosing a good basis.

 \begin{lemma}[Minkowski]\label{geom.nbres}
We can choose a basis $(Q_1, \ldots, Q_r)$ for the free part of the Mordell-Weil group satisfying
$\hat{h}(Q_1) \leq \ldots \leq \hat{h}(Q_r)$, and
\begin{equation}\label{minkowski}
\prod_{i=1}^{r}\hat{h}(Q_{i}) \leq (r!)^4 \Reg(E/K).
\end{equation}
\end{lemma}

\noindent {\it Proof. } 
Put together Minkowski's theorem on the successive minima \cite[Theorem V, Chapter VIII, section 4.3]{cassels} with Lemma 8 page 135
of \cite{cassels} as \cite[Lemma 5.1]{remond.sous2005}. 
\hfill$\Box$

\medskip
From now on, we assume that we have chosen a basis $(Q_1,\ldots,Q_r)$ as in Lemma~\ref{geom.nbres}.
Thus, in order to bound  the regulator from below, it suffices to use a lower bound for the $\hat{h}(Q_i)$'s.
In order to bound from above $\prod_{i=1}^{r}\max\{1, \hat{h}(Q_{i})\}$  (resp. $\log V = \hat{h}(Q_{r})$), we will use inequality (\ref{minkowski})
together with a lower bound for the $\hat{h}(Q_i)$'s satisfying $\hat{h}(Q_{i}) < 1$ (resp. for $\hat{h}(Q_1), \ldots,  \hat{h}(Q_{r-1})$). So we now bring our attention to the problem of giving lower bounds for the height of non-torsion points of the Mordell-Weil group.

\medskip
\noindent {\bf Lower bound for the height of non-torsion points.} 
It is known that for non-torsion points, the N\'eron-Tate height is non-zero and we can then ask for a lower bound. 
Since the elliptic curve $E/K_{0}$ is fixed, but not the degree $[K:\rat]$ of the field of rationality of the point $Q$, we are interested in Lehmer's type results. 
The following corollary of a theorem of D. Masser \cite{masser-bull-soc1989} is enough for our purpose. \footnote{Note that, as pointed out by Manin \cite{manin}, for a given elliptic curve and a given number field $K$, it is not difficult to compute an effective lower bound for $\hat{h}(Q)$.}

\begin{prop}[Masser]\label{prop-masser}
There exists a positive real number $\kappa_4 < 1$, depending on the curve $E/\Fu$,
such that, for any  field extension  $K/K_0$ of degree $d= [K:\rat] \geq 2$ and for all points $Q$ in $E(K)\setminus E_{tors}$
one has 
\begin{equation}\label{borne.masser} 
\hat{h}(Q) \geq \frac{\kappa_4}{d^{3}(\log d)^{2}}. 
\end{equation} 
\end{prop}

\medskip

We then obtain the following bounds. 

\begin{lemma}\label{bound-parameters}
The following inequalities hold :
\begin{equation}\label{bound-reg-1}
\log^+(\Reg(E/K)^{-1})\ll d^3 \, (\log^+ d)\,  (\log^+D_K) \, (\log^+\log D_K).
\end{equation}
\begin{equation}\label{bound-logV}
\log V \leq \left(\frac{d^{3}(\log^+ d)^{2}}{\kappa_4}\right)^{r-1}\cdot (r!)^4 \cdot \Reg(E/K).
\end{equation}
\begin{equation}\label{bound-prod-height-generators}
\prod_{i=1}^{r}\max\{1,\hat{h}(Q_{i})\} \leq \left(\frac{d^{3}(\log^+d)^{2}}{\kappa_4} \right)^r  \cdot(r!)^4 \cdot \Reg(E/K),
\end{equation}
where the implicit constant in the symbol $\ll$ depends at most on $E/K_0$.
\end{lemma}

\noindent {\it Proof. }
It follows from Lemma~\ref{geom.nbres}  and from Proposition~\ref{prop-masser} that we have
$$ \left(\frac{\kappa_4}{d^{3}(\log^+d)^{2}} \right)^r\leq (r!)^4\Reg(E/\Ft).$$
Hence we get, using Lemma~\ref{rank-discriminant} :
\begin{eqnarray*}
\log^+(\Reg(E/K)^{-1}) & \leq & 4 r\log r + r\log\left(\frac{d^{3}(\log^+d)^{2}}{\kappa_4}\right)
\\
& \ll & d^3 (\log^+ d)\, (\log^+D_K)\, (\log^+\log D_K).\\ 
\end{eqnarray*}
This proves (\ref{bound-reg-1}). To prove (\ref{bound-logV}), we simply
apply Masser's lower bound (\ref{borne.masser}) to the $r-1$ smallest points of the basis and replace it in Minkowski's
inequality (\ref{minkowski}). Finally, to prove (\ref{bound-prod-height-generators}),
write $\prod_{i=1}^{r}\max\{1,\hat{h}(Q_{i})\}  = \prod_{i=1}^{r}\hat{h}(Q_{i}) \times \left( \prod_{\hat{h}(Q_{i})<1}\hat{h}(Q_{i}) \right)^{-1}$,
where the second product runs over the $i$'s for which $\hat{h}(Q_{i})<1$.
Applying the inequality (\ref{minkowski}) for the first factor
and Masser's lower bound (\ref{borne.masser}) to the second one, we get the result.
\hfill$\Box$

\subsubsection{A conditional upper bound for the regulator}\label{conditional-bounds}

\medskip

\noindent {\bf On the BSD-conjecture.} 
The upper bounds (\ref{bound-logV}) and (\ref{bound-prod-height-generators}) obtained in Lemma~\ref{bound-parameters} for the height of the generators involve
the regulator $\Reg(E/K)$. In order to bound it from above, 
the BSD-conjecture suggests to bound each of the terms of the formula (\ref{formulaBSD}). 
We denote $h_{Falt}(E/K)$ the Faltings' height of $E/K$. 
The next proposition is Proposition 3.12 of \cite{mordell-weil.tate-shafa}. 
 
\begin{prop}\label{borne-regulateur} 
Suppose that Conjecture \ref{funct-eq} and \ref{bsd} hold for the elliptic curve $E/K$. 
Then  
\begin{equation} \label{borne-sha.reg-ellip}
 \Reg(E/K) \leq C_{d} \cdot D_K^{3/2} \cdot  \sqrt{N_{K/\rat}(\mathcal{F}_{E/K})}  \cdot (\exp \{ h_{Falt}(E/K)\}\cdot h_{Falt}(E/K))^{d}, 
\end{equation}
where $d = [K:\rat]$ and we may take $C_{d} = \left(\frac{9}{2\pi}\right)^{d} \cdot (3 d^2)^{d} \cdot (129. (5^d-1)(3d)^6)^{\frac{(1+3^{d/2})^8}{\log (1+3^{d/2})}}$. 
\end{prop} 

\medskip
 
Using this proposition and Lemma~\ref{conductor} for the conductor, we get a conditional bound for the regulator as we want and we can also bound
the other quantities involving the generators of the Mordell-Weil group.

\begin{lemma}\label{final-bound-generators}
Suppose that Conjecture \ref{funct-eq} and \ref{bsd} hold for the elliptic curve $E/K$. 
Then there exist positive numbers $\kappa_{5}$, $\kappa_{6}$, $\kappa_{7}$ (depending at most on $E/K_0$) such that
\begin{enumerate}
\item
$\Reg(E/K)\ll e^{\kappa_{5}^d} D_K^{3/2}.$
\item
$\log^+\log V \ll \kappa_{6}^d \,  (\log^+D_K) \, (\log^+\log D_K)$.
\item
$\displaystyle\prod_{i=1}^{r}\max\{1,\hat{h}(Q_{i})\} \ll \left( d \cdot \log ^+ D_K \right)^{\kappa_{7} d^3 \log ^+ D_K} e^{\kappa_{5}^d}\cdot D_K^{3/2},$

\end{enumerate}
where the implicit constants in the symbol $\ll$ depend at most on $E/K_0$.
\end{lemma}

\noindent {\it Proof. }
\begin{enumerate}
\item
By \cite[Remark 5.1.1, Chapter IX]{cornell-silverman}, we have $h_{Falt}(E/K) \leq h_{Falt}(E/K_0)$.
Using now Lemma \ref{conductor}, the result is an immediate consequence of Proposition~\ref{borne-regulateur}.   
\item
Use the bound (\ref{bound-logV}) of Lemma~\ref{bound-parameters}, item 1 and  Lemma~\ref{rank-discriminant}.

\item
The result follows from the bound (\ref{bound-prod-height-generators}) of Lemma~\ref{bound-parameters}, Lemma~\ref{rank-discriminant} and item 1.
\end{enumerate}
\hfill$\Box$

\medskip

\subsection{Proof of Theorem \ref{conjub}}\label{conj-u_b}

We want to apply Theorem~\ref{borne-hauteur}, so we need to estimate first $C_{E,K}$.
By Lemma~\ref{rank-discriminant}, Lemma~\ref{final-bound-generators} and Inequality (\ref{bound-reg-1}) of Lemma~\ref{bound-parameters}, we have :
\begin{eqnarray*}
\log C_{E,K}  & \leq & c_{1} d^6 \,  (\log^+ D_K)^2\, (\log (d \log^+ D_K)) + \kappa_{5}^d,
\end{eqnarray*}
for some $c_1 = c_1 (E/K_0)$. On the other hand, by Lemma~\ref{rank-discriminant} again, we have :
\begin{equation*}
8r^2+\gamma_1 dr \leq c_2 d^6 \, (\log^+D_K)^2,
\end{equation*}
for some $c_2 = c_2 (E/K_0)$. Theorem~\ref{conjub} follows at once from these estimates and from Theorem~\ref{borne-hauteur}.
\hfill $\Box$

\section{What about $abc$?}\label{appli-abc}

As already mentioned in the introduction, the second author has shown in \cite{bornes}
that one can deduce an $abc$-type inequality over number fields from a bound for the height of the $S$-integral
points on a fixed curve, explicit in the set $S$, the degree $[K:\rat]$ and the discriminant $D_K$
of the number field. The aim of this section is to prove such an inequality using the conditional bound
obtained for the integral points in Theorem~\ref{conjub},  following the method of \cite{bornes}.

With the notations of Section~\ref{notations}, the $abc$-conjecture of D. Masser \cite{masser} and J. Oesterl\'e \cite{oesterle}
over number fields can be stated  as follows (see \cite{elkies-abc}).

\begin{conj}[Masser-Oesterl\'e]\label{abc}{($abc$)} 
 
Let $\Fd$ be a number field. For every $\varepsilon >0$, there exists a 
real number $c_{\varepsilon, \Fd}>0$ such that, for  $a, b, c$ non zero
elements of $\Fd$ satisfying $a+b=c$, we have
$$h_{\Fd}(a:b:c)< (1+\varepsilon) \rad (a:b:c) +c_{\varepsilon, \Fd}.$$ 
\end{conj} 

The result that we will prove in this section is the following.

\begin{thm}\label{bsd-abc} 
 
Let $a,b,c$ be non zero elements in the number field $\Fd$ satisfying 
$a +b =c$. Let $E$ be any elliptic curve defined over some number field 
$\Fu \subset \Fd$. 
Suppose that for any finite extension $K$ of $\Fd$, the $L$-series of $E/K$ satisfies a  functional equation (Conjecture \ref{funct-eq}) 
and that the Birch and Swinnerton-Dyer Conjecture (Conjecture \ref{bsd}) holds for 
the elliptic curve $E/K$. 

Then, there exist real positive numbers $\beta_{1}$ and $\beta_{2}$ depending at most on the curve $E/\Fu$, the degree $[\Fd:\rat]$ and
the absolute value $D_{\Fd}$ of the discriminant of $\Fd$, 
such that 
$$h_{\Fd}(a:b:c)< \exp\{ \beta_{1} \rad(a:b:c)^{3} + \beta_{2} \}.$$  
 \end{thm} 
 
In Theorem \ref{bsd-abc} one may take $\beta_1 = c_1(E/\Fu) \cdot [\Fd:\rat]^6 \cdot (\log^+  [\Fd:\rat]) \cdot (\log ^+ D_{\Fd})^2$ and
$\beta_2 = c_2(E/\Fu) ^{[\Fd:\rat]}$, where $c_1(E/\Fu)$ and $c_2(E/\Fu)$ depend at most on $E/\Fu$.

Roughly speaking, the known results on $abc$ over number fields (\cite{bornes}, \cite{gyory-yu}, \cite{gyory.acta-arith2008}) give an inequality 
$$h(a:b:c) \leq \exp\{\beta_1\rad(a:b:c) + \beta_2\}.$$
In \cite{gyory.acta-arith2008}, one may take $\beta_1 = 1 + \epsilon$.
Thus the bound obtained in Theorem \ref{bsd-abc} is less good as the known results. However, our result
is obtained by a totally different method, which shows a connection between two conjectures of a very different nature,
namely the BSD-conjecture and the $abc$-conjecture. Moreover, improvements in the lower bounds used for
linear forms in elliptic logarithms could yield a better inequality, see the discussion in the remarks \ref{remark2} and
\ref{remark3} below.

\smallskip
 
Observe that the validity of the hypothesis (the functional equation and the BSD-conjecture) is only needed for a single elliptic curve $E$
(which we may choose as we wish), but for infinitely many field extensions $\Ft/\Fd$. In fact, it turns out (see Section~\ref{proof-of-bsd-abc})
that we need the hypothesis only for every extension $\Ft/\Fd$ of relative degree $[\Ft:\Fd]  \leq \deg(f)$, where $f$ is 
any fixed \Belyi \esp function associated to $E$. 
For instance, we may choose the CM curve  given by the affine equation  $y^{2} = x^{3} -x$ and for which $(x,y) \mapsto -\frac{1}{4}\frac{(1-x)^{2}}{x}$ is a \Belyi \esp map of degree 4.
Note also that if the elliptic curve has complex multiplication or is defined over $\rat$, then the 
conjecture concerning the functional equation is true 
(see \cite{shimura-taniyama} and \cite{wiles1995}). On the other hand, there is some evidence for the truth of the BSD-conjecture
(see \cite{coates-wiles1977}, \cite{gross-zagier1986}, \cite{rubin1987} and \cite{kolyvagin1988}). 

\smallskip

It is worth noting that with different methods, D. Goldfeld and L. Szpiro \cite[Theorem 2]{goldfeld-szpiro} proved that there is
a relation between the BSD-conjecture and Szpiro's conjecture. The latter one relates the discriminant of the curve with its conductor,
namely $D = O(N^{6+ \epsilon})$. It is known to imply a weak version of the $abc$-conjecture over $\rat$ (where in conjecture \ref{abc},
$1+\epsilon$ is replaced by an absolute constant).  They proved that if the order
of the Tate-Shafarevich group satisfies $|\Sha| = O(N^{1/2 + \epsilon})$, for every $\epsilon > 0$, for {\it all} elliptic curves defined over
$\rat$, then the relation $D = O(N^{18+ \epsilon})$ holds for every elliptic curve over $\rat$. 
Their proof uses the BSD-conjecture for all elliptic curves over $\rat$ in the case of rank zero, which is a theorem in this case.

\subsection{Proof of Theorem~\ref{bsd-abc}.}\label{proof-of-bsd-abc}

Let $a,b,c$ be non zero elements of the number field $\Fd$ such that $a +b =c$. 
Let $S_{1}$ be the set of the prime ideals $\vp$ of $\Fd$ such that $ \card\{\ord_{\vp}(a), \ord_{\vp}(b), \ord_{\vp}(c)\}\geq 2.$
We then have 
\begin{equation*}
\rad(a:b:c) = \sum_{\vp \in  S_{1}} \log N_{\Fd/\rat}\vp:= \Sigma_{S_1}.
\end{equation*}

To our point $(a:b:c)$ will correspond an integral point on an elliptic curve. Choose $E$ any elliptic curve defined over a subfield
$\Fu \subset \Fd$. Let
$$y^{2} = x^{3} +Ax +B$$
be a Weierstrass equation of $E$, with $A,B\in O_{\Fu}$. This curve being fixed once for all and chosen independently of $a,b,c$, all the parameters 
depending only on $\Fu$ and $E/\Fu$ will be considered as ``constants''. 
 
Using a uniformization theorem of G. V. \Belyi, we can lift the point $(a:b:c)$ to an integral point of the elliptic curve.
Indeed, by \cite[Theorem~4]{belyi}, there exists a finite surjective morphism $f : E \rightarrow \espproj^{1}$
defined over $\Fu$, unramified outside $0$, $1$ and $\infty$, and sending the origin $O = (0:1:0)$ of $E$ to $\{0, 1, \infty\}$.
Let $Q$ be a point of $E(\overline{\Fd})$ such that
$$f(Q) = (a:c) \in \mathbf{P}^{1}\setminus \{0, 1, \infty\}.$$ 
Since $b=c-a$, the point $(a:c)$ is an $S_1$-integral point of $ \espproj^{1}\setminus \{0, 1, \infty\}$ and the point $Q$ contains all the information about our triple $(a:b:c)$. Moreover, we can use the properties of the elliptic curve.

The Chevalley-Weil theorem (see \cite[\S~4.2]{serre} or \cite[Lemma 2.4 and Lemma 2.5]{bornes} for an affine version),
gives us information about the lift. 

\begin{lemma}[Chevalley-Weil]\label{chevalley-weil}
The field of definition $\Ft=\Fd(Q)$ of the point
$Q$ is a finite extension of $\Fd$ of relative degree 
\begin{equation}\label{degree}
[\Ft:\Fd] \leq \deg(f)
\end{equation}
and which is unramified outside $S = S_{1} \cup S_{0}$, for some finite set of places $S_{0}$ of $\Fd$ depending only on the curve $E/\Fu$ and
the function $f$.
Moreover, the point $Q$ is $S'$-integral, where $S'$ is the set of places of $\Ft$ lying above $S$. \footnote{See \cite{bilu-strambi-surroca}
for a quantitative version with control on the height of the set $S_0$.}
\end{lemma}

We apply now Theorem~\ref{conjub} which gives us a conditional upper bound for the height of the lift of $(a:c)$ depending on the field
extension $\Ft$ and the set of places $S'$:
\begin{equation}\label{height-pts}
h_x(Q) \leq \exp\{\alpha_1^d +\alpha_2 \,d^6 (\log^+ D_K)^2\, [\Sigma_{S'} + \log (d \log^+ D_K)],
\end{equation}
where $d=[\Ft:\rat]$ and $\alpha_1$ and $\alpha_2$ depend at most on $E/\Fu$.

To end the proof of Theorem~\ref{bsd-abc}, it remains to relate $h_{x} (Q)$ to $h(a:b:c)$ on the one hand, and the parameters
$d$, $\Sigma_{S'}$ and $D_{\Ft}$ to the radical $\rad(a:b:c)$ on the other hand. This is achieved by the following lemma. 

\smallskip

\begin{lemma}\label{inequalities}
The following inequalities hold :
\begin{enumerate}
\item
$h(a:b:c) \ll h_{x} (Q).$
\item
$d=[K:\rat] \ll [\Fd:\rat].$
\item
$\Sigma_{S'}  \ll \rad(a:b:c)$
\item
$\log D_{\Ft} \ll \rad(a:b:c) + \log D_{\Fd}.$
\end{enumerate}
where the implicit constants in the symbol $\ll$ depend at most on $K_0$ and $E/K_0$.
\end{lemma}

\noindent {\it Proof.} 
Using the basic properties of the heights and because the \Belyi \, map $f$  and the $x$-coordinate are both functions on $E$, we have 
$$h_{x} (Q) \gg h_{f}(Q) = \frac{1}{\deg(f)} h(f(Q)) = \frac{1}{\deg(f)} h(a:c) \geq \frac{1}{\deg(f)} (h(a:b:c) - \log 2).$$
Hence the first point is proved.
The second point follows from (\ref{degree}) and from the fact that the \Belyi \esp map depends only on $E/K_0$.
To prove the third item, we first note that
\begin{equation}\label{lemma-sigma}
\Sigma_S \ll \rad(a:b:c)
\end{equation}
since
$$\Sigma_{S} = \Sigma_{S_0 \cup S_1} \leq \Sigma_{S_0} + \Sigma_{S_1}=\Sigma_{S_0} + \rad(a:b:c)$$
and since $S_0$ depends only on $E/{K_0}$ (and on the choice of $f$). Now we have, by definition of $\Sigma_{S'}$:
$$\Sigma_{S'} = \sum_{\vp \in S} \sum_{\vq | \vp} \log N(\vq) \leq \sum_{\vp \in S} \sum_{\vq | \vp} \frac{f_{\vq}}{f_{\vp}} \log N(\vp)
\leq [\Ft:\Fd] \,\Sigma_{S}, $$
from which the third item follows, by (\ref{degree}) and (\ref{lemma-sigma}).
Finally, let us prove the last inequality.
From Lemma \ref{chevalley-weil}, the set of places of $\Fd$ on which the extension $\Ft/\Fd$ ramifies is contained in $S$.
We can then apply the following form of the Dedekind-Hensel inequality which is Lemma 3.17 of \cite{bilu-strambi-surroca}:
\begin{equation}
\label{discrL/K}
\log D_{\Ft} \le \Sigma_{S} + [\Ft:\Fd] \left( \log D_{\Fd} +   1.26\right).
\end{equation}
We conclude with (\ref{degree}) and (\ref{lemma-sigma}) as before.
\hfill $\Box$

\smallskip
The proof of Theorem~\ref{bsd-abc} is now straightforward : It suffices to insert the inequalities of Lemma~\ref{inequalities}
in (\ref{height-pts}).\hfill $\Box$ 

\medskip 

\subsection{Some remarks}\label{someremarks}

\begin{remark}\label{remark1}

Observe why the bound of Theorem \ref{bsd-abc} has growth order $\exp\{\rad(a:b:c)^3\}$.
The first remark is that, in the bound for the height of the integral points of the curve obtained in Theorem \ref{borne-hauteur}, the radical of $(a:b:c)$ appears in several ways. 

First, the radical appears, as expected, in the set of places $S$ (see (\ref{lemma-sigma})) :
$$\Sigma_S \ll \rad(a:b:c).$$

Next, in the bound of Theorem \ref{borne-hauteur}, the radical appears in the rank, in the height of a system of generators and in the regulator.
More precisely, by the Weak Mordell-Weil theorem, $\rk(E/K)$ can be bounded linearly by $\log D_K$ (see Lemma \ref{rank-discriminant})
which is in turn bounded linearly by the radical of $(a:b:c)$ (it comes from the method, that is, \Belyi's theorem, the Chevalley-Weil theorem and
the Dedekind-Hensel inequality; see Lemma~\ref{inequalities}, Item~4).
Thus we have 
$$\rk(E/K) \ll \log D_K \ll \rad (a:b:c).$$

Using Minkowski's theorem on the successive minima and a lower bound for the height of non-torsion points,
$\log (\Reg(E/K))^{-1}$ is bounded by $(\log D_K\,. \log \log D_K)$ (Lemma \ref{bound-parameters}), hence 
$$\log (\Reg(E/K))^{-1} \ll \rad (a:b:c)\ . \log(\rad (a:b:c)).$$
The factors concerning the heights of the generators are bounded under the BSD-conjecture by Lemma \ref{final-bound-generators}: 
$$(\log ^+ \log V)^r \leq \exp \{c_1 \rad (a:b:c) \cdot \log \rad (a:b:c)\}$$
and
$$ \prod\max \{1, \hat{h}(Q_i)\} \leq \exp \{c_2 \rad (a:b:c) \cdot \log \rad (a:b:c)\},$$
where $c_1$, $c_2$ are constants.

In the bound of Theorem~\ref{borne-hauteur} appears also a factor concerning the rank:
$$r^{2r^2} \leq \exp \{c_3 \rad (a:b:c)^2 \log \rad (a:b:c)\}.$$
Finally, we have the factors 
\begin{equation}\label{contribution-de-exp(gamma_1 d r Sigma_S)}\exp(\gamma_1 d r\Sigma_S) \leq \exp \{c_4 \rad (a:b:c)^2\}
\end{equation}
and
\begin{equation}\label{contribution-de-p^{8n(n+1)}}\exp(8 r^2 \Sigma_S) \leq \exp \{c_5 \rad (a:b:c)^3\},
\end{equation}
the latter one being the biggest contribution to  the radical. It comes from the factor $p^{8n(n+1)}$ of the bound of N. Hirata
(see Theorem~4.6 in \cite{S-entiers-Bosser-Surroca}) for the linear form in elliptic logarithms in the ultrametric case,
where $n$ is the number of logarithms (essentially $n = r$) and $p$ is the prime number lying below the ultrametric place.
\end{remark}

\medskip 

\begin{remark}\label{remark2}
According to N. Hirata, it seems possible that in the lower bound of linear forms of elliptic logarithms in the ultrametric case,
a refinement on $p$ of the order $p^n$ could be obtained (instead of $p^{8n(n+1)}$), giving a contribution of the form
$\exp \{c_6 \rad (a:b:c)^2\}$, instead of (\ref{contribution-de-p^{8n(n+1)}}). This would lead to a final bound in Theorem \ref{bsd-abc}
$$h_{\Fd}(a:b:c)< \exp\{ \beta_{1} \rad(a:b:c)^{2} \log \rad (a:b:c) + \beta_{2} \}.$$ 
The worse contribution in this case would be the factor $r^{2r^2}$ which appears both in the ultrametric and the archimedean lower bounds for
linear forms in elliptic logarithms (theorems 4.2 and 4.6 in \cite{S-entiers-Bosser-Surroca}).
It seems reasonable to conjecture that the lower bounds for linear forms in logarithms are still valid with the smaller factor $r^r$
instead of $r^{2r^2}$. In that case, the worse contributions would be the factors $p^n$ and
$\exp(\gamma_1 d r\Sigma_S)$, both yielding a factor of the shape
$\exp \{c_7 \rad (a:b:c)^2\}$. The presence of the factor $\exp(\gamma_1 d r\Sigma_S)$ shows that even a drastic improvement in the
dependence on $p$ and on the number of logarithms would not be sufficient to get a final inequality better than $\exp \{c_7 \rad (a:b:c)^2\}$.
\end{remark}

\medskip 
  
\begin{remark}\label{remark3}
S. David and N. Hirata (\cite{sinnou-noriko-Crelle2009}) suggested an elliptic analog of the classical Lang-Waldschmidt conjecture.
We quote here the following  particular case.  

\begin{conj} [Elliptic Lang-Waldschmidt] \label{lw}
Let $n$ be a rational integer $\geq 1$. There  exists a strictly positive real number $C(n)$ such that the following property holds.
Let $K$ be a number field of degree $d$ over $\rat$. Let $E/K$ be an elliptic curve given by a Weierstrass equation $y^2=x^3+Ax+B$.
Denote $h_E = \max \{1, h(1:A:B)\}$. Let $b_0, b_1, \ldots , b_n$ be rational integers and $B= \max_{1\leq i \leq n} \{ |b_i|\}$.
Let $\gamma_i$ be points in $E(K) \subset \espproj ^2(K)$ and $u_i$ be an elliptic logarithm of $\gamma_i$. Put
$\mathcal{L} = b_0 + b_1 u_1 + \cdots +  b_n u_n$. If $\mathcal{L} \not=0$, then 
$$\log | \mathcal{L}| \geq - C(n) d^2 (\log B + h_E) (\sum_{i=1}^n \max\{1, \hat h (\gamma_i) \}).$$
\end{conj}

Following the results in the classical (non elliptic) case, we can expect that $C(n) = c_1^{n}$, where $c_1$ is some
absolute constant. For the ultrametric analogue, by considering  Theorem 1' of \cite{yu-compositio}, we can conjecture a constant of the form
$C(p,n,d) = c_1^{n} \cdot p^{c_2d}$, where $c_1$ and $c_2$ are absolute. Using Conjecture~\ref{lw} with $C(n) = c_1^{n}$,
instead of David's theorem (Theorem 4.2 of \cite{S-entiers-Bosser-Surroca}), together with an ultrametric analogue of Conjecture \ref{lw}
with $C(p,n,d) = c_1^{n} \cdot p^{c_2d}$,
instead of Hirata's theorem (Theorem 4.6 of \cite{S-entiers-Bosser-Surroca}), our method would give the following $abc$-type inequality
$$h_{\Fd}(a:b:c)< \exp\{ \beta_{1} \rad(a:b:c)\cdot \log \rad(a:b:c)  + \beta_{2} \}.$$  
Indeed, with the notation of \cite{S-entiers-Bosser-Surroca}, the
factor $\exp(\gamma_1 d r\Sigma_S)$ comes from the factor $\nu_{\vp}^{2r}$, where $\nu_{\vp}$ is the exponent of a certain group satisfying
$\nu_{\vp}\leq p^{c_8d}$.
If we repeat the arguments given in \cite{S-entiers-Bosser-Surroca} for the proof of Theorem~\ref{borne-hauteur}
using Conjecture~\ref{lw}, we see that we now get $r\nu_{\vp}^2$ instead of $\nu_{\vp}^{2r}$,
and thus the factor $\exp(\gamma_1 d r\Sigma_S)$ is replaced by
$r\exp(\gamma_1 d \Sigma_S)$, giving a contribution $\exp(c_9 \rad(a:b:c))$ instead of (\ref{contribution-de-exp(gamma_1 d r Sigma_S)}).
Hence the worse contribution here would come from $\sum_{1\leq i \leq r} \max \{1, \hat{h}(Q_i)\}$ which we bound, under the BSD-conjecture, by
$$r \cdot \max_{1\leq i \leq r}\{1,\hat{h}(Q_i)\} \ll r (r!)^4 (\frac{\kappa_4}{d^3(\log d)^2})^{1-r} \cdot D_K^{3/2} \leq \exp \{c_{10}\rad (a:b:c) \cdot \log \rad (a:b:c)\}.$$
\end{remark}
   
\medskip 
   
\bibliography{mabiblio-BSD-abc}

\noindent
\textbf{Vincent Bosser}\\
Laboratoire Nicolas Oresme\\
Universit\'e de Caen\\
F14032 Caen cedex\\
France

\bigskip

 \noindent
\textbf{Andrea Surroca}\\
Mathematisches Institut\\
Universit\"at Basel\\
Rheinsprung 21\\
CH-4051 Basel\\
Switzerland

\end{document}